\documentclass[10pt]{amsart}
     \makeatletter
     \def\section{\@startsection{section}{1}%
     \z@{.7\linespacing\@plus\linespacing}{.5\linespacing}%
     {\bfseries%\normalfont\scshape
     \centering
     }}
     \def\@secnumfont{\bfseries}
     \makeatother
\newtheorem{theorem}{Theorem}[section]

\theoremstyle{definition}

\theoremstyle{remark}

\numberwithin{equation}{section}
\usepackage{amsmath,amsthm,amssymb,amsbsy}

\begin{document}

\setlength{\parindent}{0cm}
\setlength{\parskip}{0.5cm}

\title[Hitting times for jump processes]{An elementary proof that the first hitting time of an open set by a
  jump process is a stopping time}

\author{Alexander Sokol}

\address{Alexander Sokol: Institute of Mathematics, University of
  Copenhagen, 2100 Copenhagen, Denmark}
\email{alexander@math.ku.dk}
\urladdr{http://www.math.ku.dk/$\sim$alexander}



\begin{abstract}
We give a short and elementary proof that the first hitting time of an
open set by the jump process of a c\`{a}dl\`{a}g adapted process is a
stopping time.
\end{abstract}

\maketitle

\noindent

\section{Introduction}

For a stochastic process $X$ and a subset $B$ of the real numbers, the
mapping $T$ defined by $T=\inf\{t\ge0|X_t \in B\}$ is called the first hitting time
of $B$ by $X$. A classical result in the general theory of processes
is the d\'{e}but theorem, which has as a corollary that under the
usual conditions, the first hitting time of a Borel set for a progressively measurable process is a stopping time,
see \cite{DM}, Section III.44 for a proof of this theorem, or \cite{RFB1} and
\cite{RFB2} for a recent simpler proof. For many purposes,
however, the general d\'{e}but theorem is not needed, and weaker results
may suffice, where elementary methods may be used to obtain the
results. For example, it is elementary to show that the first hitting time
of an open set by a c\`{a}dl\`{a}g adapted process is a stopping time, see
\cite{PP}, Theorem I.3. Using somewhat more advanced, yet relatively
elementary methods, Lemma II.75.1 of \cite{RW1} shows that the first hitting time of a
compact set by a c\`{a}dl\`{a}g adapted process is a stopping time.

These elementary proofs show stopping time properties for the first hitting
times of a c\`{a}dl\`{a}g adapted process $X$. However, the jump process
$\varDelta X$ in general has paths with neither left limits nor right
limits, and so the previous elementary results do not apply. In this
note, we give a short and elementary proof that the first hitting time of an open set
by $\varDelta X$ is a stopping time when the filtration is
right-continuous and $X$ is c\`{a}dl\`{a}g adapted. This result may be
used to give an elementary proof that the jumps of a c\`{a}dl\`{a}g
adapted process are covered by the graphs of a countable sequence of
stopping times.

\section{A stopping time result}

Assume given a filtered probability space $(\Omega,\mathcal{F},(\mathcal{F}_t),P)$
such that the filtration $(\mathcal{F}_t)_{t\ge0}$ is right-continuous in the
sense that $\mathcal{F}_t=\cap_{s>t}\mathcal{F}_s$ for all $t\ge0$. We
use the convention that $X_{0-}=X_0$, so that there is no jump at the
timepoint zero.

\begin{theorem}
\label{theorem:FirstEntranceJumpOpen}
Let $X$ be a c\`{a}dl\`{a}g adapted process, and let $U$ be an open set in
$\mathbb{R}$. Define $T=\inf\{t\ge0|\varDelta X_t\in U\}$. Then $T$ is a
stopping time.
\end{theorem}

As $X$ has c\`{a}dl\`{a}g, $\varDelta X$ is zero everywhere except for on a
countable set, and so $T$ is identically zero if $U$ contains
zero. In this case, $T$ is trivially a stopping time. Thus, it
suffices to prove the result in the case where $U$ does not contain
zero. Therefore, assume that $U$ is an open set not containing zero. As the filtration is right-continuous, an elementary
argument yields that to show the stopping time property of $T$, it
suffices to show $(T<t)\in\mathcal{F}_t$ for $t>0$, see Theorem I.1 of
\cite{PP}.

To this end, fix $t>0$. As $X_0-X_{0-}=0$ and $U$ does not contain
zero, we have
\begin{eqnarray}
  (T<t)=(\exists\; s\in(0,\infty):s<t\textrm{ and }X_s-X_{s-}\in U)\;.
\end{eqnarray}
Let $F_m=\{x\in\mathbb{R}|\;\forall\; y\in U^c:|x-y|\ge1/m\}$, $F_m$ is an
intersection of closed sets and therefore itself closed. Clearly,
$(F_m)_{m\ge1}$ is increasing, and since $U$ is open, $U=\cup_{m=1}^\infty
F_m$. Also, $F_m\subseteq F_{m+1}^\circ$, where
$F_{m+1}^\circ$ denotes the interior of $F_{m+1}$. Let $\Theta_k$ be the subset of $\mathbb{Q}^2$ defined by $\Theta_k =
\{(p,q)\in\mathbb{Q}^2|0<p<q<t,|p-q|\le\frac{1}{k}\}$. We will prove the
result by showing that
\begin{eqnarray}
  &&(\exists\; s\in(0,\infty):s<t\textrm{ and }X_s-X_{s-}\in U)\notag\\
 &&=\cup_{m=1}^\infty\cup_{n=1}^\infty \cap_{k=n}^\infty \cup_{(p,q)\in\Theta_k} (X_q-X_p\in F_m)\;.\label{eq:Equality}
\end{eqnarray}
To obtain this, first consider the inclusion towards the right. Assume
that there is $0<s<t$ such that $X_s-X_{s-}\in U$.  Take $m$ such that
$X_s-X_{s-}\in F_m$. As $F_m\subseteq F_{m+1}^\circ$, we then have
$X_s-X_{s-}\in F_{m+1}^\circ$ as well. As $F_{m+1}^\circ$ is open and as
$X$ is c\`{a}dl\`{a}g, it holds that there is $\varepsilon>0$ such
that whenever $p,q\ge 0$ with $p\in (s-\varepsilon,s)$ and $q\in(s,s+\varepsilon)$, $X_q-X_p\in F_{m+1}^\circ$. Take $n\in\mathbb{N}$
such that $1/2n<\varepsilon$. We now claim that for $k\ge n$, there is $(p,q)\in\Theta_k$
such that $X_q-X_p\in F_{m+1}$. To prove this, let
$k\ge n$ be given. By the density properties of $\mathbb{Q}_+$ in $\mathbb{R}_+$,
there are elements $p,q\in\mathbb{Q}$ with $p,q\in(0,t)$ such that $p\in (s-1/2k,s)$ and $q\in
(s,s+1/2k)$. In particular, then $0<p<q<t$ and $|p-q|\le
|p-s|+|s-q|\le 1/k$, so
$(p,q)\in\Theta_k$. As $1/2k\le1/2n<\varepsilon$, we have $p\in (s-\varepsilon,s)$ and
$q\in (s,s+\varepsilon)$, and so $X_q-X_p\in F^\circ_{m+1}\subseteq F_{m+1}$. This proves the inclusion
towards the right.

Now consider the inclusion towards the left. Assume that there is
$m\ge1$ and $n\ge1$ such that for all $k\ge n$, there exists $(p,q)\in\Theta_k$
with $X_q-X_p\in F_m$. We may use this to obtain sequences
$(p_k)_{k\ge n}$ and $(q_k)_{k\ge n}$ with the properties that $p_k,q_k\in\mathbb{Q}$,
$0<p_k<q_k<t$, $|p_k-q_k|\le\frac{1}{k}$ and $X_{q_k}-X_{p_k}\in F_m$. Putting $p_k=p_n$ and $q_k=q_n$ for $k<n$,
we then find that the sequences $(p_k)_{k\ge1}$ and $(q_k)_{k\ge1}$
satisfy $p_k,q_k\in\mathbb{Q}$, $0<p_k<q_k<t$, $\lim_k|p_k-q_k|=0$ and
$X_{q_k}-X_{p_k}\in F_m$. As all sequences of real
numbers contain a monotone subsequence, we may by taking two
consecutive subsequences and renaming our sequences obtain the
existence of two monotone sequences $(p_k)$ and $(q_k)$ in $\mathbb{Q}$
with $0<p_k<q_k<t$, $\lim_k |p_k-q_k|=0$ and
$X_{q_k}-X_{p_k}\in F_m$. As bounded monotone sequences are
convergent, both $(p_k)$ are $(q_k)$ are
then convergent, and as $\lim_k|p_k-q_k|=0$, the limit $s\ge0$ is
the same for both sequences.

We wish to argue that $s>0$, that $X_{s-}=\lim_kX_{p_k}$ and that $X_s=\lim_k
X_{q_k}$. To this end, recall that $U$ does not contain zero, and
so as $F_m\subseteq U$, $F_m$ does not contain zero either. Also note
that as both $(p_k)$ and $(q_k)$ are monotone, the limits
$\lim_kX_{p_k}$ and $\lim_kX_{q_k}$ exist and are either equal to $X_s$ or
$X_{s-}$. As $X_{q_k}-X_{p_k}\in F_m$ and $F_m$ is closed and does not
contain zero, $\lim_kX_{q_k}-\lim_kX_{p_k}=\lim_k
X_{q_k}-X_{p_k}\neq0$. From this, we can immediately conclude that $s>0$, as if $s=0$, we
would obtain that both $\lim_k X_{q_k}$ and $\lim_kX_{p_k}$ were equal to $X_s$, yielding
$\lim_k X_{q_k}-\lim_k X_{p_k}=0$, a contradiction. Also, we cannot have that both limits are
$X_s$ or that both limits are $X_{s-}$, and so only two cases are
possible, namely that $X_s=\lim_kX_{q_k}$ and $X_{s-}=\lim_kX_{p_k}$ or that
$X_s=\lim_kX_{p_k}$ and $X_{s-}=\lim_kX_{q_k}$. We wish to argue that
the former holds. If $X_s=X_{s-}$, this is trivially the
case. Assume that $X_s\neq X_{s-}$ and that $X_s=\lim_kX_{p_k}$ and
$X_{s-}=\lim_kX_{q_k}$. If $q_k\ge s$ from a point onwards or $p_k<s$ from a point onwards,
we obtain $X_s=X_{s-}$, a contradiction. Therefore, $q_k<s$ infinitely often and $p_k\ge s$ infinitely
often. By monotonicity, $q_k<s$ and $p_k\ge s$ from a point onwards, a
contradiction with $p_k<q_k$. We conclude $X_s=\lim_kX_{q_k}$ and $X_{s-}=\lim_kX_{p_k}$, as desired.

In particular, $X_s-X_{s-}=\lim_k X_{q_k}-X_{p_k}$. As
$X_{q_k}-X_{p_k}\in F_m$ and $F_m$ is closed, we obtain $X_s-X_{s-}\in
F_m\subseteq U$. Next, note that if $s=t$, we have $p_k,q_k<s$ for all
$k$, yielding that both sequences must be increasing and $X_s=\lim
X_{q_k}=X_{s-}$, a contradiction with the fact that $X_s-X_{s-}\neq0$
as $X_s-X_{s-}\in U$.  Thus, $0<s<t$. This proves the existence of $s\in(0,\infty)$ with $s<t$ such that $X_s-X_{s-}\in U$,
and so proves the inclusion towards the right.

We have now shown (\ref{eq:Equality}). Now, as $X_s$ is
$\mathcal{F}_t$ measurable for all $0\le s\le t$, it holds that the set
$\cup_{m=1}^\infty\cup_{n=1}^\infty \cap_{k=n}^\infty
\cup_{(p,q)\in\Theta_k} (X_q-X_p\in F_m)$ is $\mathcal{F}_t$ measurable as
well. We conclude that $(T<t)\in\mathcal{F}_t$ and so $T$ is a stopping time.

\bibliographystyle{amsplain}

\begin{thebibliography}{99}

\bibitem{RFB1} Bass, Richard F.: The measurability of hitting times,
  {\em Electronic Communications in Probability} {\bf 15} (2010) 99--105.

\bibitem{RFB2} Bass, Richard F.: Correction to ``The measurability of
  hitting times'',
  {\em Electronic Communications in Probability} {\bf 15} (2010) 189-191.

\bibitem{DM} Dellacherie, C. and Meyer, P.-A.: {\em Probabilit\'{e}s et potentiel}, Hermann, 1975.

\bibitem{PP} Protter, P.: {\em Stochastic Integration and Differential
  Equations}, Springer, 2005.

\bibitem{RW1} Rogers, L. C. G. and Williams, D.: {\em Diffusions, Markov
  Processes and Martingales}, volume 1, Cambridge University Press, 2000

\end{thebibliography}

\end{document}